\documentclass[12pt]{amsart}
\usepackage{amssymb}
\newtheorem{theorem}{Theorem}[section]
\newtheorem{definition}[theorem]{Definition}

\newtheorem{lemma}{Lemma}[theorem]

\begin{document}
\commby{XXX}
\title[The maximality of the core model]{The maximality of the core model}
\author{E. Schimmerling}
\author{J.R. Steel}
\address{Department of Mathematics, University of California, 
Irvine, California 92697-3875, USA}
\email{eschimme@math.uci.edu}
\address{Department of Mathematics, University of California,
Berkeley, California 94720, USA}
\email{steel@math.berkeley.edu}
\thanks{The research of the first author 
was partially supported by a National Science Foundation
Mathematical Sciences Postdoctoral Research Fellowship at 
the Massachusetts Institute of Technology.}
\thanks{The research of the second author 
was partially supported by the National Science Foundation.}
\maketitle


\newcommand{\forces}{\Vdash}
\newcommand{\decides}{\parallel}
\newcommand\res{{{\upharpoonright}}}
\newcommand\exist{\exists}


\newcommand{\dom}{\mathrm{dom}}
\newcommand{\ran}{\mathrm{ran}}
\newcommand{\crit}{\mathrm{crit}}
\newcommand{\cf}{\mathrm{cf}}
\newcommand{\ot}{\mathrm{o.t.}}
\newcommand{\card}{\mathrm{card}}
\newcommand\otp{{\mbox{otp}}}
\newcommand{\supp}{\mathrm{supp}}
\newcommand{\support}{\mathrm{support}}
\newcommand{\lh}{\mathrm{lh}}
\newcommand{\rank}{\mathrm{rank}}
\newcommand{\id}{\mathrm{id}}
\newcommand\OR{{\mathrm{OR}}}
\newcommand\ult{{\mathrm{ult}}}
\newcommand\pred{{\mathrm{pred} \,}}
\newcommand\rt{{\mathrm{root} \,}}
\newcommand\decap{{\mathrm{decap}}}
\newcommand{\HOD}{\mathrm{HOD}}
\newcommand{\AD}{\mathrm{AD}}


\newcommand\avec{{\vec a}}
\newcommand\bvec{{\vec b}}
\newcommand\cvec{{\vec c}}
\newcommand\dvec{{\vec d}}
\newcommand\kvec{{\vec k}}
\newcommand\mvec{{\vec m}}
\newcommand\nvec{{\vec n}}
\newcommand\uvec{{\vec u}}
\newcommand\vvec{{\vec v}}
\newcommand\xvec{{\vec x}}
\newcommand\yvec{{\vec y}}

\newcommand\Avec{{\vec A}}
\newcommand\Bvec{{\vec B}}
\newcommand\Cvec{{\vec C}}
\newcommand\Dvec{{\vec D}}
\newcommand\Evec{{\vec E}}

\newcommand\gavec{{{\vec {\ga} }}}
\newcommand\gbvec{{{\vec {\gb} }}}
\newcommand\gkvec{{{\vec {\gk} }}}
\newcommand\etavec{{{\vec {\eta} }}}
\newcommand\glvec{{{\vec {\gl} }}}
\newcommand\gmvec{{{\vec {\mu} }}}
\newcommand\gnvec{{{\vec {\nu} }}}
\newcommand\gsvec{{{\vec {\gs} }}}

\newcommand\gLvec{{{\vec {\gL} }}}
\newcommand\gOvec{{{\vec {\gO} }}}

\newcommand\cMvec{{{\vec {\mathcal M} }}}
\newcommand\cNvec{{{\vec {\mathcal N} }}}
\newcommand\cPvec{{{\vec {\mathcal P} }}}
\newcommand\cQvec{{{\vec {\mathcal Q} }}}
\newcommand\cRvec{{{\vec {\mathcal R} }}}
\newcommand\cSvec{{{\vec {\mathcal S} }}}


\newcommand\htil{{\widetilde h}}
\newcommand\ktil{{\widetilde k}}

\newcommand\Ftil{{\widetilde F}}
\newcommand\Gtil{{\widetilde G}}
\newcommand\Ktil{{\widetilde K}}

\newcommand\cPtil{{\widetilde {\cP}}}
\newcommand\cQtil{{\widetilde {\cQ}}}
\newcommand\cTtil{{\widetilde {\cT}}}
\newcommand\cUtil{{\widetilde {\cU}}}
\newcommand\cVtil{{\widetilde {\cV}}}

\newcommand\gbtil{{\widetilde {\gb}}}
\newcommand\gktil{{\widetilde {\gk}}}
\newcommand\gptil{{\widetilde {\gp}}}
\newcommand\gdtil{{\widetilde {\gd}}}
\newcommand\gstil{{\widetilde {\gs}}}
\newcommand\gftil{{\widetilde {\gf}}}
\newcommand\etatil{{\widetilde {\eta}}}
\newcommand\psitil{{\widetilde {\psi}}}

\newcommand\gOtil{{\widetilde {\gO}}}


\newcommand\sdot{{\dot s}}
\newcommand\Edot{{\dot E}}
\newcommand\Fdot{{\dot F}}
\newcommand\ggdot{{\dot \gg}}
\newcommand\gmdot{{\dot \gm}}
\newcommand\gndot{{\dot \gn}}
\newcommand\gtdot{{\dot \gt}}


\newcommand\abar{{\overline a}}
\newcommand\bbar{{\overline b}}
\newcommand\cbar{{\overline c}}
\newcommand\ibar{{\overline i}}
\newcommand\pbar{{\overline p}}
\newcommand\qbar{{\overline q}}
\newcommand\rbar{{\overline r}}
\newcommand\sbar{{\overline s}}
\newcommand\tbar{{\overline t}}
\newcommand\ubar{{\overline u}}
\newcommand\vbar{{\overline v}}
\newcommand\wbar{{\overline w}}
\newcommand\xbar{{\overline x}}
\newcommand\ybar{{\overline y}}

\newcommand\Abar{{\overline A}}
\newcommand\Bbar{{\overline B}}
\newcommand\Cbar{{\overline C}}
\newcommand\Dbar{{\overline D}}
\newcommand\Ebar{{\overline E}}
\newcommand\Fbar{{\overline F}}
\newcommand\Kbar{{\overline K}}
\newcommand\Tbar{{\overline T}}
\newcommand\Wbar{{\overline W}}

\newcommand\cHbar{{\overline {\cH}}}
\newcommand\cMbar{{\overline {\cM}}}
\newcommand\cNbar{{\overline {\cN}}}
\newcommand\cPbar{{\overline {\cP}}}
\newcommand\cQbar{{\overline {\cQ}}}
\newcommand\cRbar{{\overline {\cR}}}
\newcommand\cSbar{{\overline {\cS}}}
\newcommand\cTbar{{\overline {\cT}}}
\newcommand\cUbar{{\overline {\cU}}}
\newcommand\cVbar{{\overline {\cV}}}

\newcommand\gabar{{\overline \ga}}
\newcommand\gbbar{{\overline \gb}}
\newcommand\ggbar{{\overline \gg}}
\newcommand\gdbar{{\overline \gd}}
\newcommand\gkbar{{\overline \gk}}
\newcommand\glbar{{\overline \gl}}
\newcommand\gnbar{{\overline \gn}}

\newcommand\gGbar{{\overline \gG}}
\newcommand\gObar{{\overline \gO}}
\newcommand\gUbar{{\overline \gU}}

\newcommand\cPbarbar{\overline{{\overline {\cP}}}}
\newcommand\cQbarbar{\overline{{\overline {\cQ}}}}
\newcommand\cSbarbar{\overline{{\overline {\cS}}}}
\newcommand\cTbarbar{\overline{{\overline {\cT}}}}


\newcommand\fA{\mathfrak A}
\newcommand\fB{\mathfrak B}
\newcommand\fC{\mathfrak C}
\newcommand\fM{\mathfrak M}
\newcommand\fN{\mathfrak N}
\newcommand\fR{\mathfrak R}
\newcommand\fS{\mathfrak S}
\newcommand\fT{\mathfrak T}

\newcommand\fa{\mathfrak a}
\newcommand\fb{\mathfrak b}
\newcommand\fc{\mathfrak c}
\newcommand\fm{\mathfrak m}
\newcommand\fn{\mathfrak n}
\newcommand\fr{\mathfrak r}
\newcommand\fs{\mathfrak s}
\newcommand\ft{\mathfrak t}


\newcommand{\cA}{{\mathcal A}}
\newcommand{\cB}{{\mathcal B}}
\newcommand{\cC}{{\mathcal C}}
\newcommand{\cD}{{\mathcal D}}
\newcommand{\cE}{{\mathcal E}}
\newcommand{\cF}{{\mathcal F}}
\newcommand{\cG}{{\mathcal G}}
\newcommand{\cH}{{\mathcal H}}
\newcommand{\cI}{{\mathcal I}}
\newcommand{\cJ}{{\mathcal J}}
\newcommand{\cK}{{\mathcal K}}
\newcommand{\cL}{{\mathcal L}}
\newcommand{\cM}{{\mathcal M}}
\newcommand{\cN}{{\mathcal N}}
\newcommand{\cO}{{\mathcal O}}
\newcommand{\cP}{{\mathcal P}}
\newcommand{\cQ}{{\mathcal Q}}
\newcommand{\cR}{{\mathcal R}}
\newcommand{\cS}{{\mathcal S}}
\newcommand{\cT}{{\mathcal T}}
\newcommand{\cU}{{\mathcal U}}
\newcommand{\cV}{{\mathcal V}}
\newcommand{\cW}{{\mathcal W}}
\newcommand{\cX}{{\mathcal X}}
\newcommand{\cY}{{\mathcal Y}}
\newcommand{\cZ}{{\mathcal Z}}


\newcommand\lb{\lbrack}
\newcommand\rb{\rbrack}
\newcommand\lcb{\lbrace}
\newcommand\rcb{\rbrace}


\newcommand\ra{\rightarrow}
\newcommand\la{\leftarrow}

\newcommand\thra{\twoheadrightarrow}

\newcommand\lra{\longrightarrow}
\newcommand\lla{\longleftarrow}

\newcommand\llra{\longleftrightarrow}

\newcommand\Ra{\Rightarrow}
\newcommand\La{\Leftarrow}

\newcommand\Lra{\Longrightarrow}
\newcommand\Lla{\Longleftarrow}

\newcommand\Llra{\Longleftrightarrow}
\renewcommand\iff{\Leftrightarrow}


\newcommand{\ga}{\alpha}     
\newcommand{\gb}{\beta}      
\renewcommand{\gg}{\gamma}   
\newcommand{\gd}{\delta}     
\renewcommand\ge{\varepsilon}
\newcommand{\gz}{\zeta}      
\newcommand{\gth}{\theta}    
\newcommand{\gi}{\iota}      
\newcommand{\gk}{\kappa}  
\newcommand{\gl}{\lambda}    
\newcommand{\gm}{\mu}        
\newcommand{\gn}{\nu}        
\newcommand{\gx}{\xi}        
\newcommand{\gom}{\omicron}  
\newcommand{\gp}{\pi}        
\newcommand{\gr}{\rho}       
\newcommand{\gs}{\sigma}     
\newcommand{\gt}{\tau}       
\newcommand{\gu}{\upsilon}   
\newcommand{\gph}{\varphi}      
\newcommand{\gf}{\varphi}      
\newcommand{\gch}{\chi}      
\newcommand{\gps}{\psi}      
\newcommand{\go}{\omega}     

\newcommand{\gA}{A}     
\newcommand{\gB}{B}      
\newcommand{\gG}{\Gamma}     
\newcommand{\gD}{\Delta}     
\newcommand{\gEp}{E}  
\newcommand{\gZ}{Z}      
\newcommand{\gEe}{H}      
\newcommand{\gTh}{\Theta}    
\newcommand{\gI}{I}      
\newcommand{\gK}{K}     
\newcommand{\gL}{\Lambda}    
\newcommand{\gM}{M}        
\newcommand{\gN}{N}        
\newcommand{\gX}{\Xi}        
\newcommand{\gOm}{O}  
\newcommand{\gP}{\Pi}        
\newcommand{\gR}{P}       
\newcommand{\gS}{\Sigma}     
\newcommand{\gT}{T}       
\newcommand{\gU}{\Upsilon}   
\newcommand{\gPh}{\Phi}      
\newcommand{\gCh}{X}      
\newcommand{\gPs}{\Psi}      
\newcommand{\gO}{\Omega}         

\newcommand{\bA}{{\bf A}}  
\newcommand{\bB}{{\bf B}}      
\newcommand{\bG}{\boldGamma}     
\newcommand{\bD}{\boldDelta}     
\newcommand{\bEp}{{\bf E}}  
\newcommand{\bZ}{{\bf Z}}      
\newcommand{\bEe}{{\bf H}}      
\newcommand{\bTh}{\boldTheta}    
\newcommand{\bI}{{\bf I}}     
\newcommand{\bK}{{\bf K}}     
\newcommand{\bL}{{\bf L}}    
\newcommand{\bM}{{\bf M}}        
\newcommand{\bN}{{\bf N}}        
\newcommand{\bX}{\boldXi}        
\newcommand{\bOm}{{\bf O}}  
\newcommand{\bP}{\boldPi}        
\newcommand{\bR}{{\bf P}}       
\newcommand{\bS}{\boldSigma}     
\newcommand{\bT}{{\bf T}}       
\newcommand{\bU}{\boldUpsilon}   
\newcommand{\bPh}{\boldPhi}      
\newcommand{\bCh}{{\bf X}}      
\newcommand{\bPs}{\boldPsi}      
\newcommand{\bO}{\boldOmega}     


\newcommand{\ha}{\aleph}
\newcommand{\hb}{\beth}
\newcommand{\hg}{\gimel}
\newcommand{\hd}{\daleth}
\section{Introduction}

This paper deals with the core model, or $K$, as defined in \cite{St1}. 
Our main results are interesting as part of the pure theory of this model,
but they also yield some new consistency strength lower bounds, as
we shall point out. Of course, one cannot give a sensible definition of
$K$, much less prove anything about it, without an anti-large-cardinal
hypothesis. For the most part, the hypothesis of this nature we shall
assume here is that there is no inner model satisfying "There is a
Woodin cardinal". The results we shall prove about $K$ under this hypothesis
have all been proved by others under more restrictive anti-large-cardinal
hypotheses, and our proofs build on those earlier proofs.

The term "the core model" was introduced by Dodd and Jensen to mean,
roughly, "the largest canonical inner model there is". They gave a
precise definition of the term which captures this intuition in the case
there is no inner model with a measurable cardinal in \cite{DoJ}. 
Work of Mitchell, 
and then Martin, Mitchell and Steel, led to a precise definition which
captures the intuition in the case there is no inner model with a Woodin
cardinal, and every set has a sharp.  
(Cf. \cite{St1}.  The hypothesis that every set has a sharp
seems to be a weakness in the basic theory of the core model 
for a Woodin cardinal.
Theorem~\ref{erdos} shows how a weaker Erd\"{o}s partition property suffices.)
The canonicity of the
core model is evidenced by the absoluteness of its definition: 
if $G$ is set generic over $V$, then $K^V = K^{V[G]}$.  (Cf. \cite[7.3]{St1}.)
The maximality of the core
model manifests itself in several ways. Perhaps the most important of
these is the weak covering property of $K$,
that it computes successors of singular cardinals correctly.
(Cf. \cite[8.15]{St1}, \cite{MiSchSt}, \cite{MiSch}, and \cite{SchW}.)
Other manifestations of the
maximality of $K$ are its $\Sigma^1_3$ correctness (cf. \cite[7.9]{St1}), 
and its rigidity (cf. \cite[8.8]{St1}).

Our main result here is that $K$ is maximal in another sense, namely,
that every extender which "could be added" to $K$ is already in $K$. More
precisely,
Theorem~\ref{Kmax} says that $K$ is maximal in the sense that 
any countably certified extender
that coheres with $K$ is actually on the extender sequence of $K$,
and similarly for the iterates of $K$.
Mitchell showed earlier that every countably complete
$K$-ultrafilter is in $K$ under the stronger hypothesis 
that there is no inner model with a measurable cardinal $\gk$ with
$o(\gk) = \gk^{++}$. (Cf. \cite{Mi}).
It is not known if Theorem~\ref{Kmax} holds for countably closed extenders.
Theorem~\ref{Mmax} is a result similar 
to Theorem~\ref{Kmax} for mice with sufficiently
small projectum.

Our first application of Theorem~\ref{Kmax} is Theorem~\ref{wc}, 
which says that $K$ computes successors of
weakly compact cardinals correctly.
This extends Kunen's result for $L$ 
and the weak covering properties of $K$ mentioned above.

Recall from \cite{St1} 
that $K$ is the transitive collapse of an elementary substructure of $K^c$,
the background certified core model, and that $K^c$ is the limit of
the maximal, $1$-small construction $\mathbb{C}$.
In Theorem~\ref{KtoKc}, another application of Theorem~\ref{Kmax},
we show that $K^c$ and all the models on $\mathbb{C}$ are iterates of $K$.

Theorem~\ref{KtoKc} is used to prove Theorem~\ref{mahlo},
which says, roughly, that a Mahlo cardinal is a kind of
closure point of the construction $\mathbb{C}$
(although not in as strong a sense as we would like).

The universality of initial segments of $K$ is related to some open questions
in descriptive set theory, such as the consistency strength of $u_2 = \ha_2$.
Jensen showed that it is consistent for a countable mouse to out-iterate
$K \| \ha_1$.
Theorem~\ref{omega-2}, another application of maximality,
shows that $K \| \gk$ is universal for mice of
height $\gk$ whenever $\gk$ is a cardinal $\geq \ha_2$.
(We proved Theorem~\ref{omega-2} assuming countable closure,
and Mitchell saw how to eliminate the assumption using \cite{MiSch}.)

Using Theorem~\ref{wc} and \cite{Sch2}, we prove Theorem~\ref{swc},
which says that if $\gk$ is a weakly compact cardinal
and $\square^{<\go}_\gk$ fails, 
then there is an inner model with a 
Woodin cardinal.  
The definition of $\square^{<\go}_\gk$
is recalled in Section~\ref{section-square};
it is weaker than Jensen's principle $\square_\gk$,
and stronger than his $\square^*_\gk$.
Woodin noted that our proof,
when combined with his core model induction,
improves the conclusion of Theorem~\ref{swc}
to $L(\mathbb{R})$-determinacy.  

The proof of Theorem~\ref{swc}, 
together with the results in \cite{Sch2} and \cite{MiSchSt},
are used in the proof of Theorem~\ref{sssl}, where we show that
if $\gk$ is a singular, strong limit cardinal, then there is an inner model
with a Woodin cardinal.  With little additional work, the conclusion can be
strengthened to the existence of an inner model with Woodin cardinals
cofinal in its ordinals, which, by Woodin, implies Inductive Determinacy.

Some of the theorems in this paper
have the hypothesis that $\gO$ is a measurable cardinal.
This is really just being used to see that 
the theory of \cite{St1} applies
(sharps for bounded subsets of $\gO$ would do).
For our theorems that have this hypothesis on $\gO$,
Theorem~\ref{erdos} shows how to reduce it
to the partition property:
$\gO \lra (\go )^{< \go }_\alpha$ for all $\ga < \gO$.
For example, if it exists, the Erd\"{o}s cardinal 
$\gk ( \go )$ has this property,
and it is consistent with $V = L$ that $\gk ( \go )$ exists.

$K$ is a particular model of the form $L[\Evec ]$ where $\Evec$ is a coherent 
sequence of extenders.  
Section 6 pertains to models of the form $L[\Evec ]$ in general.
With no anti-large cardinal hypothesis,
we show that if $L[\Evec ]$ is adequately iterable, 
then it is also maximal in that $L[\Evec]$ satisfies the statement,
"if $F$ is a countably certified extender which coheres with $\Evec$,
then $F = E_\ga$ for some $\ga$".
One effect of this is that, typically, if $L[\Evec ]$ satisfies a large
cardinal property, then the large cardinal property is witnessed by
extenders on $\Evec$.
\section{Maximality}

\begin{definition}\label{certificate}
Suppose that $\cM$ is a premouse and that $F$ is an
extender that coheres with $\cM$.
Let $\gk = \crit (F)$ and $\gn = \gn (F)$.
Suppose that 
$$\cA \subseteq \bigcup_{n<\go} \cP ( [ \gk ]^n )^\cM .$$
Then a 
{\bf weak $\cA$-certificate for $(\cM , F )$}
is a pair $(N,G)$ such that 
\begin{list}{}{}
\item[(a)]{$N$ is a transitive, power admissible set, $^\go N \subset N$,
$$\gk \cup \cA \cup \{ \cJ^\cM_\gk \} \subset N,$$ 
and $G$ is an extender over $N$,}
\item[(b)]{$F \cap \left( [\gn ]^{<\go} \times \cA \right)
= G \cap \left( [\gn ]^{<\go} \times \cA \right),$}
\item[(c)]{$V_{\gn +1} \subset \ult (N , G)$,}
\item[(d)]{if $j : N \lra \ult ( N,G)$ is the ultrapower map, then $\cM$
and $j(\cJ^\cM_\gk)$ agree below $\lh (F)$, and}
\item[(e)]{$\langle \cP (\ga ) \cap N , \in \rangle \prec
\langle \cP (\ga ) , \in \rangle$ whenever $\ga < \gk$.}
\end{list}
\end{definition}

Weak $\cA$-certificates
are weaker than the $\cA$-certificates of \cite[1.1]{St1} 
in two ways.  
First, we weakened the condition that $V_\gk \subseteq N$
to just condition (e).  And, second, 
$\ult (N,G)$ is not required to be countably closed.

The following definition should be compared with \cite[1.2]{St1}.

\begin{definition}\label{certified}
Suppose that $\cM$ is a premouse and $F$ is an extender of length $\ga$
that coheres with
$\cM$.  Then $(\cM , F)$ is {\bf weakly countably certified} 
iff for all countable
$$\cA \subseteq \bigcup_{n<\go} \cP ( [ \gk ]^n ) \cap J^\cM_\ga$$
there is a weak $\cA$-certificate for $(\cJ^\cM_\ga , F)$.
\end{definition}

When $\cT$ is an iteration tree, we shall sometimes write
$\cM (\cT , \eta )$ for $\cM^\cT_\eta$, the $\eta$'th model of $\cT$,
to avoid double superscripts.
When $\cT$ has a final model, we shall denote it 
either by $\cM^\cT_\infty$ or by $\cM (\cT , \infty )$.
The following is our main maximality theorem.

\begin{theorem}\label{Kmax}
Suppose that $\gO$ is a measurable cardinal and that there is no inner model
with a Woodin cardinal.  
Let $\cT$ be a normal
iteration tree on $K$ of successor length $< \gO$.  
Suppose that $F$ is an extender that coheres with
$\cM^\cT_\infty$ and that $$\lh (E^\cT_\eta ) < \lh (F)$$ for all 
$\eta < \lh ( \cT)$.
Suppose that $(\cM^\cT_\infty , F)$ is weakly countably certified.  
Then $F$ is on the $\cM^\cT_\infty$-sequence.
\end{theorem}

The following notion shall enter into the proof of Theorem~\ref{Kmax},
and is of independent interest.

\begin{definition}\label{conditionally-strong}
A premouse $\cM$ is {\bf conditionally $\ga$-strong} iff
there is a normal iteration tree $\cT$ on an $\ga$-strong premouse
such that $\cT$ has successor length,
and there is an elementary embedding $\gp$ of $\cM$ into
an initial segment of $\cM^\cT_\infty$ with $\crit ( \gp ) \geq \ga$.
\end{definition}

Note that by a L\"{o}wenheim-Skolem argument,
both $\cT$ and $\cM$ can be taken to have the same cardinality as $\cP$
in Definition~\ref{conditionally-strong}.  

\begin{proof}[Proof of Theorem~\ref{Kmax}]
Fix $\gO$, $\cT$ and $F$ as in the statement of Theorem~\ref{Kmax}
and let $\gk = \crit (F)$, $\gn = \gn ( F)$, and $\gl = \lh (F)$.

Let $\gm_0$ be an inaccessible cardinal $< \gO$ 
such that $F \in V_{\gm_0}$.
Let $W$ be an $A_0$-soundness witness for $\cJ^K_{\gm_0}$.
Then $\cT$ induces an iteration tree on
$W$ with the same extenders, drops, degrees, and tree structure as $\cT$,
and it is enough to show that $F$ is on the sequence of the final model of
the induced tree.  For the rest of the proof, 
we shall write $\cT$ for the induced iteration tree on $W$, 
and never again refer to the original iteration tree on $K$.

As in \cite[9.7]{St1},
we write $\Phi ( \cT )$ for the phalanx derived from $\cT$.
The $F$-extension of $\Phi ( \cT )$ is defined in \cite[8.5]{St1}.
By \cite[8.6]{St1}, we are done if we show that the
$F$-extension of $\Phi ( \cT)$ is $(\gO + 1)$-iterable.
(Here we use the hypothesis that $\cT$ is normal.)
\ For contradiction, suppose that $\cU$ is an illbehaved iteration tree on the
$F$-extension of $\Phi ( \cT )$.  We interpret $\cU = \emptyset$
to mean that the ultrapower in the definition of the $F$-extension of
$\Phi (\cT )$ is illfounded.
By the tree property at $\gO$, $\lh ( \cU ) < \gO$.

By the argument of \cite[6.14]{St1},
we may fix a successor cardinal $\gm$ such that $\gm_0 < \gm < \gO$ and 
\begin{itemize}
\item{there is an iteration tree $\cS$ on 
$\cJ^W_\gm$ with the same
extenders, drops, degrees, and tree structure as $\cT$, such that
$\cM^{\cS}_\eta$
is an initial segment of $\cM^\cT_\eta$ whenever $\eta
< \lh ( \cT )$, and}
\item{there is an illbehaved iteration tree $\cV$
on the $F$-extension of $\Phi (\cS )$
with the same 
extenders, drops, degrees, and tree structure as $\cU$, 
such that $\cM^{\cV}_\eta$ is an initial
segment of $\cM^\cU_\eta$ whenever $\eta < \lh ( \cU )$.}
\end{itemize}

Recall that the $F$-extension of $\Phi (\cS)$ is the phalanx 
$$\Phi ( \cS ) ^\frown 
\left\langle \ult_k (\cP ,F), k ,\gn  \right\rangle $$
where $\cP$ is the initial segment of $\cM^\cS_\gg$ to which $F$ would be
applied according to the rules for $\go$-maximal iteration trees, and 
$k < \go$ is the degree.
(This is not literally true, 
as the official definition of a phalanx, \cite[9.5]{St1}, 
has a fourth coordinate, in our case, either $\gn$ or $\gl$,
depending on the type of $\cP$; we shall suppress this fourth coordinate.)

Let $X$ be a countable elementary substructure of $V_\gO$ with 
$\cV \in X$ and let
$$\gp_X : M_X \lra V_\gO$$ be the corresponding uncollapse.
Say $\gp_X ( \cS_X ) = \cS$, $\gp_X ( F_X ) = F$, 
$\gp_X ( \gn_X ) = \gn$, 
$\gp_X ( \gl_X ) = \gl$, 
$\gp_X ( \cP_X ) = \cP$, and $\gp_X ( \cV_X ) = \cV$.
By the absoluteness argument like that in the proof of \cite[6.14]{St1},
$\cV_X$ is a simple, illbehaved iteration tree on
$$\Phi ( \cS_X ) ^\frown 
\left\langle \ult_k (\cP_X ,F_X ),k,\gn_X \right\rangle ,$$
the $F_X$-extension of $\Phi ( \cS_X )$.

Let $Y$ be an elementary substructure of $V_\gO$ with 
$$X \cup (\gl +1) \subset Y$$ and $\card (Y ) = \card ( \gn )$.
Let $\gp_Y : M_Y \lra V_\gO$ be the uncollapse of $Y$ and 
$$\gs = \gp_Y^{-1} \circ \gp_X.$$
Say $\gp_Y \left( \cS_Y \right) = \cS$ and 
$\gp_Y \left( \cP_Y \right) = \cP$.
Note that $\cS_Y$ has the same extenders, drops, degrees, and 
tree structure as $\cS$, and that $\cM (\cS_Y , 0)$ 
is $\crit (\gp_Y )$-strong.

Let $$\cA = \ran ( \gs ) \cap \bigcup_{n<\go} \cP ( [ \gk ]^n )^{\cP_Y}$$ and 
$(N,G)$ be a weak $\cA$-certificate for $(\cP_Y ,F)$.
This is possible because $(\cM^\cT_\infty , F)$ is weakly countably
certified, and $\cP_Y$ and $\cP$ agree with
$\cM^\cT_\infty$ on subsets of $\gk$.
Let $j : N \lra \ult (N,G)$ be the ultrapower map.

Pick a coordinate $b \in [\gn ]^{< \go}$ together with functions 
$u \mapsto \gn_{u \ },$
$u \mapsto \gl_{u \ },$
and 
$u \mapsto \gs_{u \ }$
in $N$ such that
$$\gn = [b , u \mapsto \gn_u ]^N_G \ ,$$
$$\gl = [b , u \mapsto \gl_u ]^N_G \ ,$$
and
$$\gs = [b , u \mapsto \gs_u ]^N_G \ .$$
This is possible by clause (c) in Definition~\ref{certificate}.
By {\L}o{\'{s}}' theorem,
there is a set $B \in G_b$ such that for all $u \in B$,
\begin{itemize}
\item{$\gs_u \res |\cM |$ is an elementary embedding of $\cM$ into 
$\gs_u (\cM )$ for every initial segment 
$\cM$ of a model on $\Phi (\cS_X )$,}
\item{$\gn_u = \gs_u ( \gn_X ) < \gk$, and}
\item{$\gl_u = \gs_u ( \gl_X  ) < \gk$.}
\end{itemize}

\begin{lemma}\label{Kmax1}
For $G_b$-almost every $u \in B$,
there is a $\Sigma_{k+1}$-elementary, weak $k$-embedding
$$\gt : \ult_k (\cP_X , F_X ) \lra \cP_Y$$
such that $\gt \res \gl_X  = \gs_u \res \gl_{X \ }$.
\end{lemma}

\begin{proof}
Let $\langle f_n \mid n < \go \rangle$ be a sequence of 
$\gS_k (\cP_Y )$-functions
and $$\langle b_n \mid n < \go \rangle$$ be a sequence of coordinates from
$[\gl ]^{<\go}$ such that
$$\left| \ult_k ( \cP_Y , F ) \right| \cap \ran ( \gs )
= \left\{ \left[ b_n , f_n \right]^{\cP_Y}_F \mid n < \go \right\}$$
We assume that $b  = b_0 \subset b_n \subset b_{n+1}$ for all $n < \go$.
We assume that $f_n$ is chosen to be a constant function whenever possible.
We assume that if $$\left[ b_n \ , f_n \right]^{\cP_Y}_F 
= \gs ( \xi )$$ for some $\xi < \gl_X  $,
then $\gs ( \xi ) \in b_n$ and $f_n$ is the projection function
$u \mapsto u^{\gs ( \xi ), a }$.

For each $\xi < \gl_X $ and $n \in \go$ with $\gs (\xi ) \in b_n$,
let 
$$C^n_\xi = \left\{ u \in [\gk ]^{|b_n|}  \mid
\gs_{u^{b_0 , b_n}} \left( \xi \right) 
= u^{\gs(\xi ) , b_n} \right\} .$$
Then $C^n_\xi \in G_{b_n}$.
Since $N$ is countably closed (or, just because
the map $u \mapsto \gs_u$ and all the reals are in $N$),
it follows that
$$\left\langle C^n_\xi \mid \xi < \gl_X   \ \land \ \gs ( \xi ) \in b_n 
\right\rangle$$
is an element of $N$.

If $n < \go$ and $e$ is the G\"{o}del number for a $\gS_k$ formula 
$\varphi$ such that 
$$\ult_k ( \cP_Y , F ) \models \varphi 
\left[ \left[b_n \ , f_n \right]^{\cP_Y}_F \right],$$
then let
$$T^n_e
= \left\{ 
u \in [\gk ]^{|b_n |} 
\mid \cP_Y \models \varphi [ f_n (u) ]  \right\} .$$
Then, 
$T^n_e \in F_{b_n} \cap \cA$
by {\L}o\'{s}' theorem and our choice of $\cA$.
Moreover, the partial function $(n,e) \mapsto T^n_e$
is an element of $N$ since $N$ is countably closed.

Let $\left\langle  B_n  \mid n < \go \right\rangle$ be a sequence in $N$
such that $B_n \in G_{b_n}$ for all $n < \go$,
with the following three properties.
First, $B_0 =B$.
Second, for every $m, e < \go$ such that $T^m_e$ is defined,
there is an $n \geq m$ with
$$( B_n )^{b_m , b_n} \subseteq T^m_e \ .$$
Third, for every $m < \go $ and $\xi < \gl_X $,
there is an $n \geq m$ such that
$$(B_n)^{b_m , b_n} \subseteq C^m_\xi \ .$$

We need a ``fiber'' through the sets $B_n$.
There is a tree 
$\mathbb A$ in $N$ that searches for a function $g: \go \lra [\gk]^{<\go}$
such that 
$$g(m) \subset g(n) \in B_n$$
and
$$(b_n )^{g(m),g(n)} = b_m$$
whenever $m < n < \go$.
Since $\langle b_n \mid n < \go \rangle$ is a branch through 
$j ( {\mathbb A} )_b$, there is a branch through $j ( {\mathbb A} )_b$ in
$\ult (N , G)$.  Hence, for $G_b$-almost every $u \in B$,
there is a branch $g$ through ${\mathbb A}_u$ in $N$.

Given $u \in B$, fix $g \in N$ such that $g$ is a branch
through ${\mathbb A}_u$, and
define a map
$$\gt : \ult_k ( \cP_X , F_X ) \lra \cP_Y$$ by
$$\gt \left( \gs^{-1} \left( \left[ b_n , f_n \right]^{\cP_Y}_F \right) 
\right)
= f_n \left( g \left( n \right) \right).$$
Because we ``met'' the sets $T^n_{e \ }$,
$\gt$ is a weak $k$-embedding of $\ult_k ( \cP_X  , F_X )$ into $\cP_Y$
(recall that weak $k$-embeddings need not be $\gS_{k+1}$-elementary).
Because $f_n$ was chosen to be the constant function whenever possible,
if $i_X : \cP_X \lra \ult_k ( \cP_X , F_X )$ is the ultrapower map,
then $\gt \circ i_X = \gs$.  Since $i$ is a $k$-embedding and $\gs$ is
elementary, $\gt$ is a $\gS_{k+1}$-elementary.
Because we ``met'' the sets 
$C^n_{\xi \ }$, it follows that
$$\gt \res \gl_X  = \gs_u \res \gl_X .$$
\end{proof}

\begin{lemma}\label{Kmax1.5}
There are functions $u \mapsto \cP_u$ and $u \mapsto \gt_u$ in $N$
such that for $G_b$-almost every $u \in B$, both $V$ and $N$ satisfy
\begin{list}{}{}
\item[(a)]{$\cP_u$ is a conditionally $(\gl_u + 1)$-strong premouse,}
\item[(b)]{$\cP_u$ and $\cJ_\gk^\cP$ below $\gl_u +1$, and}
\item[(c)]{$\gt_u$ is a $\gS_{k+1}$-elementary, weak $k$-embedding of
$$\ult_k ( \cP_X , F_X )$$
into $\cP_u$ such that $\gt_u \res \gl_X = \gs_u \res \gl_X$.}
\end{list}
\end{lemma}

\begin{proof}
First note that $\cP_Y$ is conditionally $(\gl + 1 )$-strong as witnessed
by an initial segment of the iteration tree $\cS_Y$ and the identity map
on $\cP_Y$.

Given $u \in B$ and a corresponding $\gt$ as in Lemma~\ref{Kmax1},
let
$$\cP^*_u = \cH^{\cP_Y}_\go 
\left( (\gl_u +1 ) \cup \ran \left( \gt \right) \right),$$
and let $$\gt^*_u : \ult_k ( \cP_X  , F_X ) \lra \cP^*_u$$ be the 
corresponding collapse of $\gt$.
Then the functions 
$u \mapsto \cP^*_u$ and  $u \mapsto \gt^*_u$
satisfy conditions (a), (b), and (c) in $V$,
but these functions are not necessarily in $N$.

By the elementarity given in clause (e) of  Definition~\ref{certificate},
for each appropriate $u \in B$,
there is a premouse $\cP_u$ and an embedding $\gt_u$,
both in $N$, which satisfy
conditions (a), (b), and (c) of Lemma~\ref{Kmax1.5} in both $N$ and $V$.  
Since $N$ satisfies AC,
there are appropriate functions 
$u \mapsto \cP_u$ and $u \mapsto \gt_u$ in $N$.
\end{proof}

Fix functions as in Lemma~\ref{Kmax1.5}
and let
$$\cP^* = [b , u \mapsto \cP_u ]^N_G$$
and 
$$\gt^* = [b, u \mapsto \gt_u ]^N_G \ .$$

\begin{lemma}\label{Kmax2}
$\gt^*$ is a $\gS_{k+1}$-elementary, weak $k$-embedding of 
$$\ult_k ( \cP_X  , F_X )$$ into $\cP^*$,
and $\gt^* \res \gl_X  = \gs \res \gl_X $.
\end{lemma}

\begin{proof}
Immediate from the definitions of $\cP^*$ and $\gt$,
clause (c) of Lemma~\ref{Kmax1.5}, and {\L}o{\'{s}}' theorem.
\end{proof}

\begin{lemma}\label{Kmax3}
$\cM ( \cS_Y , \infty )$ and $\cP^*$ agree below $\gl$.
\end{lemma}

\begin{proof}
Clearly $\cM ( \cS_Y , \infty )$, $\cM(\cS , \infty)$, and
$\cM( \cT ,\infty )$ all agree beyond $\gl$.
By coherence, $\cM( \cS_Y , \infty )$
and $\ult ( \cP_Y , F )$ agree below $\gl$.
Clearly $\cM (\cS_Y , \infty )$ and $\cP_Y$ agree below $\gk$.
Hence, by clause (d) of Definition~\ref{certificate},
$\cM ( \cS_Y ,  \infty )$ and $j ( \cJ^{\cP_Y}_\gk )$ agree below $\gl$.
By clause (b) of Lemma~\ref{Kmax1.5},
$\cP_u$ agrees with $\cJ^\cP_\gk$ beyond $\gl_u$
for $G_b$-almost every $u \in B$.
It follows by {\L}o{\'{s}}' theorem that
$\cM ( \cS_Y , \infty )$ agrees with $\cP^*$ below $\gl$.
\end{proof}

\begin{lemma}\label{Kmax4}
$\cP^*$ is conditionally $(\gl +1)$-strong.
Moreover, if $\cS^*$ is an iteration tree which
witnesses that $\cP^*$ is conditionally $(\gl + 1)$-strong,
then $\cS^*$ uses the same extenders as $\cS_Y$ below $\gl$.
\end{lemma}

\begin{proof}
By {\L}o{\'{s}}' theorem and clause (a) of Lemma~\ref{Kmax1.5},
$\cP^*$ is conditionally $(\gl + 1)$-strong in $\ult (N , G)$.
By clause (c) of
Definition~\ref{certificate}, $V_{\gn +1} \subset \ult (N,G)$.
It follows that $\cP^*$ is conditionally $(\gl + 1)$-strong in $V$.
In light of Lemma~\ref{Kmax3},
any iteration tree $\cS^*$ witnessing
that $\cP^*$ is conditionally $(\gl + 1)$-strong must
use the same extenders as $\cS_Y$ below $\gl$.
\end{proof}

\begin{lemma}\label{Kmax5}
$\Phi ( \cS_Y ) ^\frown \langle \cP^* , k , \gn \rangle$ is an iterable
phalanx.
\end{lemma}

\begin{proof}
Put $\fA = \Phi ( \cS_Y ) ^\frown \langle \cP^* , k , \gn \rangle$.
That $\fA$ is a phalanx follows from Lemma~\ref{Kmax3}.

Fix $\cS^*$ and an embedding $\psi$ of $\cP^*$ into
an initial segment 
of $\cM( \cS^* , \infty )$ which witness that $\cP^*$ is conditionally
$(\gl + 1)$-strong.

Let us write $\psi ( \cP^* )$ for the initial segment of 
$\cM( \cS^* , \infty)$ into which $\psi$ elementarily embeds $\cP^*$.
Since $\cM ( \cS^* , 0 )$ is $(\gl + 1)$-strong, there is
\begin{itemize}
\item{a successor-length iteration tree $\cR$ on $W$,
such that $\lh (E^\cR_\eta ) > \gl$ whenever $\eta < \lh ( \cR )$, and}
\item{an elementary embedding $\varphi$ of
$\cM ( \cS^* , 0 ) $ into an initial segment of $\cM^\cR_\infty$
such that $\crit ( \varphi ) > \gl$.}
\end{itemize}
We may use $\varphi$ to copy $\cS^*$ to an iteration tree $\varphi \cS^*$ 
on $\cM^\cR_\infty$.  Let 
$$\varphi_\infty : \cM ( \cS^* , \infty )  \lra 
\cM ( \varphi \cS^* , \infty )$$ 
be the final map in the aforementioned copying construction.
Then $\crit ( \varphi_\infty ) > \gl$ and $\varphi \cS^* $ uses the same
extenders as $\cS_Y$ below $\gl$.
Therefore
$\varphi_\infty \circ \psi$ is an elementary embedding of $\cP^*$ into
the initial segment $(\varphi_\infty \circ \psi ) ( \cP^* )$ of 
$\cM ( \varphi \cS^* , \infty )$ and
$$\crit (\varphi_\infty \circ \psi ) > \gl .$$

Recall that $\cT$ was our original iteration tree on $W$.
Put $$\fB = 
\Phi ( \cT ) ^\frown \langle ( \varphi_\infty \circ  \psi )( \cP^* ) 
, k , \gn \rangle$$
The agreement described above implies that $\fB$ is a phalanx.
It is a phalanx ``derived'' from 
$\cT$ and $\cR ^\frown \varphi \cS^*$,
both iteration trees on $W$.
Such $W$-generated phalanxes are iterable by \cite[6.9]{St1}.

But now, the pair of maps 
$(\gp_Y , (\varphi_\infty \circ \psi ) )$
can be used
in the standard way to reduce the iterability of $\fA$ to that of $\fB$.
\end{proof}

Recall that $\cV_X $ is a simple, illbehaved iteration tree on
$$\Phi ( \cS_X  ) ^\frown \langle \ult_k ( \cP_X  , F_X ) , k , \gn_X  
\rangle$$ 
Using the pair of maps $(\gs  , \gt^* )$ we can copy $\cV_X $ to an illbehaved
iteration tree on 
$\Phi ( \cS_Y ) ^\frown \langle \cP^* , k , \gn \rangle$.
But now we have a contradiction of Lemma~\ref{Kmax5}.  This completes the
proof of Theorem~\ref{Kmax}.
\end{proof}

We conclude this section with a version of Theorem~\ref{Kmax} for mice
with sufficiently small projectum.
Note that, unlike in Theorem~\ref{Kmax}, we do not assume that 
there is no inner model with a Woodin cardinal,
nor any other smallness condition, in Theorem~\ref{Mmax}.

\begin{theorem}\label{Mmax}
Let $\cM$ be a sound premouse such that
$$\gr_\go ( \cM ) \leq \crit (E)$$
whenever $E$ is an extender from the $\cM$-sequence.
Suppose that every countable phalanx which is realizable in an $\cM$-based 
phalanx is $\go_1$-iterable.
Let $\cT$ be a normal iteration tree on $\cM$.  Suppose that $F$ is 
an extender that coheres with $\cM^\cT_\infty$ such that 
$$\lh (E^\cT_\eta ) < \lh (F)$$
for all $\eta < \lh ( \cT )$, and such that 
$$\gr_\go ( \cM ) < \crit (F).$$
Suppose that $(\cM^\cT_\infty , F)$ is weakly countably certified.
Then $F$ is on the $\cM^\cT_\infty$-sequence.
\end{theorem}

\begin{proof}
Put $\gk = \crit (F)$, $\gn = \gn (F)$, and $\gl = \lh (F)$.
Let $\theta$ be a large regular cardinal and
$\gp : N \lra V_\theta$
be the uncollapse of a countable elementary substructure of $V_\theta$.
Say $\gp ( \cTbar ) = \cT$, $\gp (\Fbar ) = F $, $\gp ( \gnbar ) = \gn$,
and $\gp ( \glbar ) = \gl$.
By the argument in the proof of Theorem~\ref{Kmax},
there is an iteration tree $\cT^*$ on $\cM$ such that 
$\cM ( \cT^* , \infty )$ agrees with $\cM^\cT_\ga$ below $\gn$ 
whenever $\ga < \lh ( \cT )$
and a map $\gt$ from 
$\cM ( \cTbar , \infty )$ into $\cM ( \cT^* , \infty )$
such that $\gt \res \gl = \gp \res \gl$.
Moreover,
if $k$ is the degree of the $F$-extension of $\Phi ( \cT )$, 
then $\gt$ is a weak $k$-embedding which is $\gS_{k+1}$-elementary.
This implies that the $\Fbar$-extension of $\Phi ( \cTbar )$ is 
realizable in the $\cM$-based phalanx
$$\Phi ( \cT ) ^\frown \left\langle \cM ( \cT^* , \infty ) , k , \gn
\right\rangle .$$
By the iterability hypothesis on $\cM$,
there is a successful coiteration $(\cU , \cV )$
of $\Phi ( \cTbar )$ and the $\Fbar$-extension of $\Phi ( \cTbar )$.
Just as in the proof of \cite[8.6]{St1},
our hypothesis on the projectum of $\cM$ allows us to conclude 
that $\Fbar$ is on the $\cM ( \cTbar , \infty )$-sequence.
(Note that we do not know that $(\cU , \cV ) \in N$ until
after arguing that $\Fbar$ is on the $\cM ( \cTbar , \infty )$-sequence.)
By the elementarity of $\gp$, $F$ is on the $\cM^\cT_\infty$-sequence.
\end{proof}
\section{Applications of maximality}

Kunen showed that if $\gk$ is a weakly compact cardinal
and $$(\gk^+ )^L < \gk^+ ,$$ then $0^\#$ exists; see \cite[p.384]{J}.
Our first application of Theorem~\ref{Kmax}
is an extension of Kunen's result to the current generation of core models.

\begin{theorem}\label{wc}
Suppose that $\gO$ is a measurable cardinal, and that there is no inner
model with a Woodin cardinal.  Let $\gk$ be a weakly compact cardinal
$<\gO$.  Then $(\gk^+)^K =\gk^+$.
\end{theorem}

The result analogous to Theorem~\ref{wc} was proved 
for measurable cardinals $\gk$ in
\cite[8.15]{St1}, for countably closed, singular cardinals in 
\cite{MiSchSt}, and for all singular cardinals in \cite{MiSch}.  
From \cite{MiSch}, we already
know that $(\gk^+)^K$ has cofinality at least $\gk$
whenever $\gk$ is a cardinal $\geq \ha_2$.

\begin{proof}[Proof of Theorem~\ref{wc}]

We shall write $K^A$ for $K^{L[A]}$ for any set $A$.
Our first few lemmas will be used in later sections and 
do not require that $\gk$ be a weakly compact cardinal.

\begin{lemma}\label{wc1}
Let $\gk$ be a cardinal.
Suppose that $A$ is a bounded subset of $\gO$
such that $(H_\gk)^{L[A]} = H_{\gk \ }$.
Then $K^A$ and $K$ agree below $\gk$.  
Moreover, for any $K$-cardinal $\ga < \gk$ 
and any properly small premouse $\cP \in L[A]$:
$$\mbox{$\cP$ is $\ga$-strong $\Llra$
($\cP$ is $\ga$-strong)$^{L[A]}$.}$$
\end{lemma}

For the definitions, see \cite[6.4, 6.12]{St1}.
To prove Lemma~\ref{wc1},
one must note that $\gO$ is an $L[A]$-indiscernible, 
and that all the theorems of \cite{St1} hold for $K$ 
constructed inside $L[A]$.
Then Lemma~\ref{wc1}
is immediate from the following minor strengthening of
\cite[6.14]{St1}.

\begin{lemma}\label{wc2}
Suppose that $\ga$ is a cardinal of $K$.
Let $\cP$ be a properly small premouse such that 
$\cJ_\ga^\cP = \cJ_\ga^K$,
but $\cP$ is not $\ga$-strong.
Then there exists a properly small premouse $\cQ$ 
which is $\gb$-strong for
every $\gb < \ga$, such that
$$\card ( | \cQ | ) = \card(\ga ),$$
and there exists a countable, simple, illbehaved iteration tree
$\cT$ on $((\cQ , \cP ), \ga)$.
\end{lemma}

The only difference between Lemma~\ref{wc2} and \cite[6.14]{St1} 
is that we do not assume that $\card(|\cP |) = \card (\ga )$.
But the proof of \cite[6.14]{St1} gives this with no additional work.
Since Lemma~\ref{wc2} holds, not only in $V$, but also in any $L[A]$,
Lemma~\ref{wc1} follows.

\begin{lemma}\label{wc3}
Suppose that $\gk$ is a cardinal of uncountable cofinality, and let
$\cP$ be a properly small premouse with 
$\cJ_\gk^\cP  = \cJ^K_{\gk \ }$.
Then: $$\mbox{$\cP$ is $\gk$-strong $\Llra$ 
$\forall$ $K$-cardinal $\ga < \gk$ ($\cP$ is $\ga$-strong).}$$
\end{lemma}

\begin{proof}
The $\Lra$ direction is obvious from the definition; 
see \cite[6.4]{St1}.
So suppose that $\cP$ is not $\gk$-strong.
By Lemma~\ref{wc2}, there is an $\gk$-strong, properly small
premouse $\cQ$ of cardinality $\gk$ such that 
$\cJ_\gk^\cQ  = \cJ^K_\gk$,
and a countable, illbehaved iteration tree $\cT$ on
$((\cQ , \cP ),\gk)$.
Since $\cT$ is countable and $\gk$ has uncountable cofinality,
there is an $\ga < \gk$ such that
$\cT$ can be construed as an iteration tree on the phalanx
$((\cQ , \cP ),\ga)$.  But then $\cP$ is not $\ga$-strong by
\cite[6.11]{St1}.
\end{proof}

Also, Lemma~\ref{wc3} holds in $L[A]$ for any set $A \in V_\gO$.
Directly from \cite[6.4]{St1}, 
we see that both $V$ and $L[A]$ satisfy the sentence:
``Every properly small, $\gk$-strong premouse $\cP$ with 
$\cJ_\gk^\cP  = \cJ^K_\gk$ and $\gr_\go (\cP ) = \gk$ is a level of $K$,
and for cofinally many $\gg < ( \gk^+ )^K$,
$\cJ^K_\gg$ is a properly small, $\gk$-strong premouse which projects to
$\gk$.''
But then, it follows that:

\begin{lemma}\label{wc4}
Let $\gk$ be a cardinal of uncountable cofinality.
Suppose that $A$ is a bounded subset of $\gO$
such that $(H_\gk)^{L[A]} = H_\gk$.
Then $K^A$ and $K$ agree below $(\gk^+)^{K^A}$.
Moreover, for any properly small premouse $\cP \in L[A]$:
$$\mbox{$\cP$ is $\gk$-strong $\Llra$
($\cP$ is $\gk$-strong)$^{L[A]}$.}$$
\end{lemma}

The hypothesis of uncountable cofinality is not needed in
Lemma~\ref{wc4}, as we shall show in Lemma~\ref{squarelemma}.

Now fix a weakly compact cardinal $\gk < \gO$
and put $\gl = (\gk^+ )^K$. Assume, for contradiction, that $\gl < \gk^+$.

Let $N$ be a transitive, power admissible set such that 
$^\go N \subset N$, 
$$V_\gk \cup \cJ^K_{\gl+1} \subset N ,$$
and $\card (N) = \gk$.
Let $A \subseteq \gk$ such that $N \in L[A]$.

It follows immediately from Lemma~\ref{wc4} that
$\gl = \left( \gk^+ \right)^{K^A}$, and that $K^A$ and $K$ agree below $\gl$.

Since $A^\#$ exists, $\left( \gk^+ \right)^{L[A]} < \gk^+$,
so $\card \left( \cP (\gk ) \cap L[A] \right) = \gk$.
Because $\gk$ is weakly compact, 
there is a non-principal, $\gk$-complete 
$L[A]$-ultrafilter $U$ over $\gk$; see \cite[p.384]{J}.
Then $\ult (L[A],U)$ is well-founded.
Let $$j: L[A] \lra \ult(L[A], U) = L[j(A)]$$ be the ultrapower map.
Then $\crit (j) = \gk$ and $A = j(A) \cap \gk \in L[j(A)]$,
so $L[A] \subseteq L[j(A)]$.
We remark that $\cP (\gk ) \cap L[j(A)] \not\subset L[A]$ is possible;
see \cite[p.394]{J}.
It follows from Lemma~\ref{wc4}
that $\gl = (\gk^+ )^{K^{j(A)}}$, and that
$K^{j(A)}$ and $K$ agree past $\gl$.
Thus $\cP (\gk ) \cap K^{j(A)} = \cP (\gk ) \cap K^A$.

Let $E_j$ be the superstrong extender derived from $j$.
So, for every $a \in [j(\gk )]^{<\go}$,
$$(E_j )_a 
= \left\{ x \subseteq [\gk ]^{|a|} \mid 
x \in L[A] \ \land \ a \in j(x) \right\} .$$
We claim that for every sequence 
$\langle \ x_\ga \mid \ga < \gk \ \rangle$ 
in $L[A]$,
$$E_j \cap \left( 
[j(\gk )]^{<\go} \times 
\left\{ x_\ga \mid \ga < \gk \right\} 
\right) 
\in L[j(A)].$$
The argument is due to Kunen, as in \cite[1.1]{MiSt}:
just note that for any $a \in [j(\gk )]^{<\go}$,
$$E_a = \left\{ x_\ga \subseteq [\gk ]^{|a|} \mid
\ga < \gk \ \land \ a \in j(x_\ga ) \right\}$$
and 
$$\left\langle j (x_\ga ) \mid \ga < \gk \right\rangle \in L[j(A)] .$$

Set 
$$F = E_j \cap \left( [j(\gk )]^{<\go} \times K^A \right)$$
and 
$$G =  E_j \cap \left( [j(\gk )]^{<\go} \times N \right).$$
Then $(K^{j(A)}, F)$ and $(N , G)$ are elements of $L[j(A)]$.
Moreover, $L[j(A)]$ satisfies the sentence,
``$F \res \xi$ coheres with $K$ and 
$(N,G)$ is an $\cA$-certificate for 
$(K , F \res \xi )$ where
$$\cA = \bigcup_{n<\go} \cP ([\gk ]^n )^K$$
and unbounded many $\xi < j(\gk )$''.
(For clause (c) of Definition~\ref{certificate},
simply note that $L[A] \models$``$V_\gk \subset N$'',
so 
$$L[j(A)] \models\mbox{``}V_{j(\gk)} \subset j(N)\mbox{'',}$$
and $j(N)$ and $\ult (N, G)$ have the same sets of rank $< j(\gk )$.)
That $\gO$ is an $L[j(\gk )]$-indiscernible is enough to see that the
conclusion of Theorem~\ref{Kmax} holds in $L[j(\gk )]$.  Therefore
every initial segment of $F$ is on the $K^{j(A)}$ sequence.
But then $\gk$ is a Shelah cardinal in $K^{j(A)}$, contradicting that there
is no inner model with a Woodin cardinal.
\end{proof}

By definition,
$K$ is the transitive collapse of an elementary substructure of
$K^c$; see \cite[\S5]{St1}.
Our next application of 
Theorem~\ref{Kmax} shows that there is an
iteration tree on $K$ with last model $K^c$.
For the statement of Theorem~\ref{KtoKc}, recall
that $K^c$ is defined as the limit of the maximal $1$-small construction
$\langle \cN_\xi \mid \xi < \gO \rangle$; see \cite[\S1]{St1}.

\begin{theorem}\label{KtoKc}
Assume that $\gO$ is a measurable cardinal and that there is no inner model
with a Woodin cardinal.
Let $\ga \leq \gO$ and $k < \go$.
Suppose that $(\cU,\cV)$ is the coiteration of 
$$\left( K , \fC_k ( \cN_\xi ) \right).$$
Then $\cV$ is trivial.
In particular, since $K$ is universal, $K^c$ is an iterate of $K$.
\end{theorem}

\begin{proof}
The proof is by induction.  Fix $\xi \leq \gO$ and $k < \go$.
Assume that $\fC_\ell ( \cN_\zeta )$ does not move in its coiteration with
$K$ whenever $$(\zeta , \ell ) <_{lex} (\xi , k ).$$
For contradiction, assume that $\fC_k (\cN_\xi )$ moves
in its coiteration with $K$.

Let $(\cU , \cV)$ be the coiteration of 
$$\left( K , \fC_k ( \cN_\xi ) \right)$$
and $\theta < \lh ( \cV )$ be least such that 
$E^\cV_\theta \not= \emptyset$.
Put $\cT = \cU \res (\theta +1)$ and $F = E^\cV_\theta$.
Then $F$ coheres with $\cM^\cT_\theta$ but $F$ is not on the
$\cM^\cT_\theta$-sequence.

First suppose that $k = 0$
and that $F$ is the last extender of $\cN_\xi$.
Then $F$ is countably certified and $\crit (F)$ is an inaccessible cardinal,
by the definition of $\cN_\xi$.
Hence, by Theorem~\ref{KtoKc}, $F$ is on the $\cM^\cT_\theta$-sequence, a
contradiction.

Next, suppose that $k = 0$ but that $F$ is not the last extender of
$\cN_\xi$.  Then there is some $\zeta < \xi$ and $\ell < \go$ 
such that $F$ is on the $\fC_\ell ( \cN_\zeta )$-sequence.  
From the fact that $\fC_\ell ( \cN_\zeta )$ and $\fC_k ( \cN_\xi )$ 
agree below $\lh (F)+1$, it follows that 
$\fC_\ell ( \cN_\zeta )$ moves in its coiteration with $K$.  
This contradicts our induction hypothesis.

So $k > 0$.  By the induction hypothesis, 
$$\fC_{k-1} (\cN_\xi ) \not= \fC_k ( \cN_\xi ).$$
Thus, $\fC_{k-1} ( \cN_\xi )$ is $(k-1)$-sound but not $k$-sound.
By the induction hypothesis,
there is an iteration tree $\cS$ on $K$ such that the last model of $\cS$ is 
$\fC_{k-1} ( \cN_\xi )$.
But then $\fC_k ( \cN_\xi )$ is a model along the main branch of $\cS$.
Namely,
$$\fC_k ( \cN_\xi ) = (\cM_{\gg +1}^*)^\cS$$ 
where $\gg +1$ is the last drop to an ultrapower of degree at most $k-1$
along the main branch of $\cS$.
But then
$\cS \res (\gg + 1)$ witnesses that the theorem holds for $(\xi , k)$,
a contradiction.
\end{proof}

Recall that $\cN_\gk$ is a premouse of height $\gk$ whenever $\gk$ 
is a cardinal, but that $\cN_\gk$ need not equal $\cJ^{K^c}_\gk$.
Theorem~\ref{mahlo} appears to be the first step in showing that 
$\cN_\gk = \cJ^{K^c}_\gk$ when $\gk$ is a Mahlo cardinal, 
but we do not know if this is true.

\begin{theorem}\label{mahlo}
Suppose that there is no inner model with a Woodin cardinal and that $\gO$
is a measurable cardinal.  Let $\kappa < \gO$ be a Mahlo cardinal.
Then there is an iteration tree $\cT$ on $\cJ^K_\gk$ such that 
$\cM^\cT_\infty = \cN_\gk$.  In particular, 
$$(\ga^+ )^{\cN_\gk} = \ga^+$$
for all singular cardinals $\ga < \gk$,
so $\cN_\gk$ is universal for mice of height $\leq \gk$.
\end{theorem}

\begin{proof}
By Theorem~\ref{KtoKc}, there is an
iteration tree $\cT$ on $K$ with $\cM^\cT_\infty = \cN_\gk$.
Certainly $\lh ( \cT ) \leq \gk + 1$.
Let us assume that $\lh ( \cT ) = \gk + 1$, the other case being clear.
It is enough if we show that $\cT$ does not drop and that
$i^\cT_{0,\gk} ( \ga ) < \gk$ for all $\ga < \gk$.

Suppose otherwise.
Form an internally approachable elementary chain 
$$\langle X_\ga \mid \ga < \gk \rangle$$ of submodels of 
$\langle V_{\gk+1} , \in , \cT \rangle$ 
such that $\card (X_\ga ) = \ga$ and 
$^\go X_\ga \in X_{\ga +1}$ for all $\ga < \gk$.
Let $$\gp_\ga : P_\ga \lra V_{\gk +1}$$ be the uncollapse of $X_\ga$ and 
$G_\ga$ be the length $\gk$ extender derived from $\gp_\ga$ 
whenever $\ga < \gk$.
There is a stationary set of inaccessible cardinals $\ga < \gk$ such that
if $(\gb +1 ) T \gk$ and $\ga = \pred^\cT (\gb +1 )$,
then
$\ga = \crit ( \gp_\ga )$,
$\gp_\ga ( \ga ) = \gk$,
$X_\ga$ is countably closed,
$\cM^\cT_\ga \cap \cP( \ga ) \in P_\ga$,
and $(P_\ga , G_\ga )$ is an $\cA$-certificate for
$(\cM^\cT_\gb , E^\cT_\gb )$ where 
$$\cA = \bigcup_{n<\go}  \cP ([\ga ]^n  ) \cap | \cM^\cT_\gb | .$$
The proof is similar to the standard argument that coiterations of mice
terminate.  Fix any such $\ga < \gk$.
Clearly $E^\cT_\gb$ coheres with $\cN_\gk$ and 
$(\cN_\gk , E^\cT_\gb )$ is countably certified.
Note that $\lh ( E^\cT_\gb )$ is a cardinal in $\cN_{\gk \ }$.
Let $\gg < \gk$ be the supremum of the $\eta < \gk$ such that 
$\gr_\go ( \cN_\eta ) < \lh ( E^\cT_\gb )$.  It is not hard to see
that $\cN_\gg$ is the passive initial segment of $\cN_\gk$ of height 
$\lh ( E^\cT_\gb )$,
$\cN_\gg {\ }^\frown \langle E^\cT_\gb \rangle$ is a premouse,
and that $(\cN_\gg , E^\cT_\gb )$ is certified as above.
Therefore
$$\gr_\go ( \cN_{\gg +1} ) \leq \gn (E^\cT_\gb ) < \lh ( E^\cT_\gb ),$$
a contradiction.
\end{proof}

\begin{theorem}[with Mitchell]\label{omega-2}
Suppose that there is no inner model with a Woodin cardinal and that 
$\gO$ is a measurable cardinal.  Let $\gk$ be a cardinal such 
that $\ha_2 \leq \gk < \gO$.
Then $\cJ^K_\gk$ is universal for premice of height $\gk$.
\end{theorem}

Theorem~\ref{omega-2} 
reduces quickly to the case in which $\gk$ is a successor
cardinal.  The authors proved Theorem~\ref{omega-2} under the
assumption that $\gk$ is a countably closed, successor cardinal.
Then, Mitchell saw how to eliminate the countable closure,
by adapting the methods of \cite{MiSch}.

\begin{proof}[Sketch of Theorem~\ref{omega-2}]
Suppose that $\cM$ is a counterexample to the theorem
and that $(\cS , \cT)$ is the coiteration of $(\cJ^K_\gk , \cM)$.
Then, there is club set $C \subseteq [0,\gk]_S$ such that
$$\cD^\cT \cap [\ga , \gk]_T = \emptyset ,$$
$$\crit( i^\cT_{\ga , \gk } ) = \ga ,$$
and $$i^\cT_{\ga , \gk} ( \ga ) = \gk$$
whenever $\ga \in C$.

Obtain $(P_\ga , G_\ga )$ from an internally approachable elementary chain
of length $\gk$ as in the proof of Theorem~\ref{mahlo}.
Here, $\card (P_\ga ) = \gm$ and it is possible that no
$P_\ga$ is countably closed.

Assuming that $\gm^\go < \gk$, it is possible to arrange countable closure.
Then, as in the proof of Theorem~\ref{mahlo},
for a stationary subset of $C$, $(P_\ga , G_\ga)$ is a weak 
$\cA$-certificate for $(\cM^\cT_\ga , E^\cT_\ga )$ where 
$$\cA = \bigcup_{n<\go}  \cP ([\ga ]^n  ) \cap | \cM^\cT_\gb | .$$
It follows that $(\cM^\cS_\ga . E^\cT_\ga )$ is weakly countably certified.
(Note that $(P_\ga , G_\ga )$ cannot be an $\cA$-certificate
since $V_\ga \not\subseteq P_\ga$.) 
\ By Theorem~\ref{Kmax}, $E^\cT_\ga$ is on the $\cM^\cS_\ga$-sequence,
an obvious contradiction.

Now suppose that $\gm^\go \geq  \gk$.
The only difference is that, with the same definitions as above,
we cannot arrange that $P_\ga$ is countably closed.
But notice that in the proof of Theorem~\ref{Kmax},
countable closure was really only used
to see that the $E^\cT_\ga$-extension of $\Phi ( \cS )$ is iterable
(so that we could apply \cite[8.6]{St1}).  
The method of \cite{MiSch} adapts to show
that for a stationary set of $\ga < \gk$,
the $E^\cT_\ga$-extension of $\Phi ( \cS )$ is iterable.
\end{proof}
\section{Square principles}\label{section-square}

Recall from \cite[5.1]{Sch1} or \cite{Sch2} the following hierarchy of
principles intermediate to 
Jensen's well-known
principles $\square_\gk$ and $\square^*_{\gk \ }$ from \cite{J1}.
If $1 < \gl \leq \gk^+$, 
then $\square_\gk^{< \gl}$ is the principle 
asserting the existence of a sequence
$\langle \cF_\gn \mid \gn  < \gk^+ \rangle$
such that for every limit ordinal $\gn \in (\gk , \gk^+)$,
\begin{list}{}{}
\item[(1)]{$1 \leq \card ( \cF_\gn ) < \gl$, and}
\item[(2)]{for all $C \in \cF_\gn$,
\begin{list}{}{}
\item[(a)]{$C$ is a closed, unbounded subset of $\gn$,}
\item[(b)]{$C$ has order type $\leq \gk$, and}
\item[(c)]{$C \cap \gm \in \cF_\gm$ whenever $\gm < \gn$
is a limit point of $C$.}
\end{list}}
\end{list}
We write $\square_\gk^\gl$ for $\square_\gk^{<\gl^+}$ ;
so $\square_\gk \equiv \square_\gk^1$ and $\square^*_\gk \equiv
\square_{\gk \ }^\gk$.  
Using these forms of weak square, 
some lower bounds on the large cardinal consistency
strength of the Proper Forcing Axiom 
and other principles were obtained in 
\cite{Sch1, Sch2, SchSt}.

Theorem~\ref{wc} 
and its proof will be used to prove Theorem~\ref{swc} below.

\begin{theorem}\label{swc}
Suppose that there is no inner model with a Woodin cardinal.
Then $\square_\gk^{<\go}$ holds whenever $\gk$ is a 
weakly compact cardinal.
\end{theorem}

If there is a measurable cardinal $\gO$, 
but no inner model with a Woodin
cardinal, then $\square_\gk^{<\go}$ holds for 
every weakly compact cardinal
$\gk < \gO$; this is immediate from 
Theorem~\ref{wc} and \cite{Sch2}.  
Theorem~\ref{swc} says that the measurable cardinal is not necessary.

Our proof of Theorem~\ref{swc} can be modified to show that if 
$\square_\gk^{<\go}$
fails, then every $x \subseteq \gk$ is in a proper class inner model 
$M_1 (x)$ with a Woodin cardinal, 
and that $M_1 (x)^\#$ exists.  It follows
from unpublished work of Woodin that $\neg \square_\gk^{<\go}$ implies 
{\boldmath $\Sigma$}$^1_2$-determinacy; 
see \cite{N} for related results.
Using our
proof of Theorem~\ref{swc}, Woodin has since shown that 
$\neg \square_\gk^{<\go}$
implies $L(\mathbb R )$-determinacy.
Woodin's proof is an induction on the levels of $L(\mathbb R )$, with 
Theorem~\ref{swc} 
as the base case, and he makes the general case look like
the base case.  Woodin had used this method before, 
in unpublished work.

\begin{proof}[Proof of Theorem~\ref{swc}]
Suppose that $\gk$ is a weakly compact cardinal and
that there is no inner model with a Woodin cardinal.
Assume for contradiction that $\square_\gk^{< \go}$ fails.
Then $A^\#$ exists whenever $A$ is a bounded subset of $\gk$.
This follows from Jensen's theorem that $\square_\gk$ holds
in $L[A]$ (cf. \cite{J1}) and Kunen's theorem
on successors of weakly compact cardinals relativized to $L[A]$
(cf. \cite[p.384]{J}).
Since $\gk$ is weakly compact,
$A^\#$ exists for all $A \subseteq \gk$.
Let $\gO = \gk^+$.
Then $\gO$ is an $A$-indiscernible for every $A \subseteq \gk$.  So, for
every $A \subseteq \gk$, the theory of $K$ up to $\gO$ of \cite{St1} applies
in $L[A]$.

Let $\left\langle A_\ga \mid \ga < \gO \right\rangle$
be a sequence of subsets of $\gk$ and
$\left\langle \gl_\ga \mid \ga < \gO \right\rangle$
be a sequence of ordinals with $A_\ga \in L[A_\gb ]$,
$$\left( \gk^+ \right)^{L[A_\ga ]} = \gl_\ga \geq \ga \ ,$$
and $\gl_\ga < \gl_\gb$ whenever $\ga < \gb < \gO$.
Also, assume that $V_\gk^{L[A_\ga ]} = V_{\gk \ }$.

Note that $\gk$ is weakly compact in $L[A_\ga ]$ whenever $\ga < \gO$.
Therefore,
$$\gl_\ga = (\gk^+ )^{K^{A_\ga}}$$
whenever $\ga < \gO$,
by a version of Theorem~\ref{wc} applied inside $L[A_\ga ]$.
By Lemma~\ref{wc4},
$K^{A_\ga}$ and 
$K^{A_\gb}$ agree below $\gl_\ga$ whenever $\ga < \gb < \gO$.
Let
$$\Evec = \bigcup_{\ga < \gO} \left( \Edot^{K^{A_\ga}} \res \gl_\ga
\right)$$
and 
$\cM = \cJ^\Evec_\gO$.
Clearly, $$(\gk^+ )^\cM = \gO .$$
By \cite{Sch2}, 
$\square_\gk^{< \go}$ holds in $K^{A_\ga}$ for all $\ga < \gO$.
The proof in \cite{Sch2} is ``local'' and gives more.
Namely, for each $\ga < \gO$,
there a sequence 
${\cF^\ga} 
= \left\langle \cF^\ga_\gn \mid \gn  < \gl_\ga \right\rangle$
that witnesses $\square_\gk^{< \go}$ in $K^{A_\ga}$,
such that $\cF^\ga = \cF^\gb \res \gl_\ga$
whenever $\ga < \gb < \gO$.
So $\cF = \bigcup_{\ga < \gO}  \cF^\ga$
witnesses $\square_\gk^{< \go}$ in $\cM$, and therefore in $V$.
This is a contradiction.
\end{proof}

If there is a measurable cardinal $\gO$, but no inner model with a 
Woodin cardinal, then $\square_\gk^{<\go}$ holds for every singular cardinal 
$\gk < \gO$.  This is immediate from \cite{MiSch} and \cite{Sch2}.
The next theorem says that the measurable cardinal is not necessary, at
least when $\gk$ is a strong limit cardinal.

\begin{theorem}\label{sssl}
Suppose that there is no inner model with a Woodin cardinal.
Then $\square_\gk^{<\go}$ holds whenever $\gk$ is a singular, 
strong limit cardinal.
\end{theorem}

Jensen proved that if $0^{\P}$ does not exist, then
$\square_\gk$ holds whenever $\gk$ is a singular cardinal.
(Cf. \cite{J2} and \cite{J3}.)
It is not known if Theorem~\ref{sssl} holds for all singular cardinals,
nor if its conclusion can be strengthened to $\square_\gk$ holding.
Woodin's inductive proof of determinacy mentioned before adapts to
show that $\neg \square_\gk^{<\go}$ implies Inductive Determinacy,
hence Projective Determinacy.
(Using the lower-part closure $Lp$ instead of $L$ in the proof, 
one can see directly that 
if $\square_\gk^{<\go}$ fails and $\gk$ is a singular, strong limit cardinal,
then, over any bounded subset of $\gk$,
there is a premouse $\cM$ with $\omega$-many Woodin cardinals 
cofinal in $\OR^\cM$.)
It is not known if $\neg \square_\gk^{<\go}$ implies
$L(\mathbb R)$-determinacy.

\begin{proof}[Proof of Theorem~\ref{sssl}]
Suppose that $\gk$ is a singular, strong limit cardinal and
that there is no inner model with a Woodin cardinal.
Assume for contradiction that $\square_\gk^{< \go}$ fails.
Jensen showed that $\square_\gk$ holds in $L[A]$ 
whenever $A \subseteq \gk$ (cf. \cite{J1})
Therefore $(\gk^+ )^{L[A]} < \gk^+$ whenever $A \subseteq \gk$.

\begin{lemma}\label{sharps}
$A^\#$ exists whenever $A \subseteq \gk$.
\end{lemma}

\begin{proof}[Sketch]
It is clear that Jensen's covering lemma 
(cf. \cite{DJ}) can be adapted to show that
$A^\#$ exists whenever $A$ is a bounded subset of $\gk$.
Suppose that $A$ is an unbounded subset of $\gk$
and that $\gl = (\kappa^+ )^{L[A]} < \kappa^+$.
We shall get $A^\#$ by piecing together lift-ups of sharps for bounded 
subsets of $\kappa$.

Let $\pi :  N \lra V_\theta$ for some large $\theta$
with $N$ transitive, countably closed, $\card (N)<\kappa$,
and $\pi$ cofinal in $\gl$.
Say $\pi( \Abar ) = A$, $\gp ( \gkbar ) = \gk$, and $\gp (\glbar ) = \gl$.
Let $E$ be the length $\kappa$ extender derived from $\pi$.
If $\glbar$ is not a cardinal of $L[\Abar ]$, 
then as usual we let $\cM$ be the first level of
$L[\Abar ]$ which sees this, note that 
$$\cM = \cH^\cM_{n+1} ( \gkbar \cup p_{n+1}^\cM )$$
for an appropriate $n < \go$ (fine structure works above $\gkbar$), 
set $\cM^* = \ult_n( \cM, E)$, and get
that $\cM^*$ is a level of $L[A]$ collapsing $\lambda$.

So $\glbar$ is a cardinal of $L[\Abar ]$, and thus $E$ measures all
subsets of $\gkbar$ in $L[\Abar]$.  
Let $\pi^*: L[\Abar ] \lra \ult ( L[\Abar ],E) =
L[A]$ be the lift-up of $\pi$.
Let $D$ be the class of strong limit cards of cofinality $\kappa^+$. 
Note the tuples from $D$ are
$\Abar$-indiscernible with respect 
to parameters $< \gkbar$. Also, $\pi^* ( \alpha )=\alpha$ for
$\alpha \in D$.
By the elementarity of $\pi^*$, the tuples from $D$ are $A$-indiscernible
with respect to parameters from $\kappa$. Since $D$ is stationary, we
are done.
\end{proof}

Let $\gkvec = \langle \gk_i \mid i < \cf ( \gk ) \rangle$ be an increasing 
sequence of cardinals cofinal in $\gk$.
Put $\gO = \gk^+$.
Choose $\left\langle A_\ga \mid \ga < \gO \right\rangle$ and 
$\left\langle \gl_\ga \mid \ga < \gO \right\rangle$ as in the proof of 
Theorem~\ref{swc},
with the additional stipulation that $\gkvec \in L[A_0]$.
Then $\gk$ is a singular, strong limit cardinal in 
$L[A_\ga ]$, so
$$\gl_\ga = ( \gk^+ )^{K^{A_\ga}}$$
for all $\ga < \gO$,
by a version of \cite{MiSchSt} applied inside $L[A_\ga ]$.
If $\gk$ has uncountable cofinality, then we are done as in the proof
of Theorem~\ref{swc}.  So assume that $\cf ( \gk ) = \go$.
Then Lemma~\ref{squarelemma} below is clearly enough.

\begin{lemma}\label{squarelemma}
If $\ga < \gb < \gO$, then 
$K^{A_\ga}$ and $K^{A_\gb}$ agree below $\gl_\ga$.
\end{lemma}

\begin{proof}
Suppose, for contradiction, that the lemma is false.
Throughout the proof, we work in $L[A_\gb]$.  Put $A = A_\ga$.

Since $K^A$ and $K$ agree below $\gk$,
there is a properly small level $\cP$ of $K^A$ 
such that $\cP$ is not $\gk$-strong.

Let $\cQ$ be a properly small premouse of cardinality $\gk$ such that
$((\cQ , \cP ) , \gk)$ is not iterable,
and let $\cT$ be an illbehaved iteration tree on $((\cQ , \cP ) , \gk)$.
By taking elementary substructures,
we find a sequence of properly small premice
$$\cQvec = \langle \cQ_i \mid i < i < \go \rangle$$
such that $\cQ_i$ is $\gk_i$-strong and 
$\card \left(\left| \cQ_i \right| \right) = \gk_i$ whenever $i < \go$,
but the phalanx
$(( \cQvec , \cP ), \gkvec )$
is not iterable.
Each $\cQ_i$ elementarily embeds into $\cQ$ with critical point $> \gk_i$.
Moreover, $\cT$ can be construed as an iteration tree on 
$((\cQvec , \cP ), \gkvec )$, and thus construed, $\cT$ is illbehaved.
Let $M$ be the transitive collapse of a countable elementary substructure of
$V_\gO$ and let
$\gp : M \lra V_\gO$ be the uncollapse.  Assume that 
$$\gp \left( ( (\cQvec^M , \cP^M ),  \gkvec^M ) \right) 
= ((\cQvec , \cP ), \gkvec )$$
and $\gp ( \cT^M ) = \cT$.
By elementarity and absoluteness,
$\cT^M$ is an illbehaved iteration tree on $((\cQvec^M , \cP^M ), \gkvec^M )$.

Now $\cQ_i \in L[A]$ for $i < \go$
since $V_\gk \subset L[A]$.
But there is no reason to believe that $\cQvec \in L[A]$.
By \cite[6.14]{St1}, $\cQ_i$ is $\gk_i$-strong in $L[A]$
whenever $i < \go$.
Because $L[A]$ contains all the reals, 
$M \cup \{ M \} \subset L[A]$.

There is a tree in $L[A]$ that simultaneously searches for:
\begin{itemize}
\item{a sequence $\cRvec$ such that $((\cRvec, \cP ), \gkvec )$
is a phalanx of proper premice and $\cR_i$ is $\gk_i$-strong in $L[A ]$
whenever $i< \go$, and}
\item{a sequence $\gsvec$ of length $\go +1$ such that
$\gs_\go$ is an elementary embedding of $\cP^M$ into $\cP$ with
$\gs_\go ( \gk^M ) = \gk$,
and $\gs_i$ is an elementary embedding of $\cQ_i^M$ into $\cR_i$
with $\gs_i ( \gk_i^M ) = \gk_i$
whenever $i < \go$,
such that
$$\gs_i \res (\gk^M_i  +1 ) = \gs_j  \res (\gk^M_i +1 )$$
whenever $i < j \leq \go$.}
\end{itemize}
There is an infinite branch through this tree
determined by $\cRvec = \cQvec$ and
$$\gsvec = \left\langle \ \gp \res \left| \cQ_i^M \right| 
\mid i < \go \right\rangle
{\ }^\frown \ 
\left\langle \gp \res \left| \cP^M \right| \right\rangle.$$
By the absoluteness, there is a branch through this tree in $L[A ]$.
So let $\cRvec$ and $\gsvec$ be determined by such a branch.

The iteration tree on $((\cRvec , \cP ), \gkvec )$ obtained by copying 
$\cT^M$ using $\gsvec$ we denote by 
$\gsvec ''\cT^M$.  Now $((\cRvec , \cP ), \gkvec )$ is iterable in
$L[A]$ since $\cP$ and the $\cR_i$ are appropriately strong in
$L[A]$.
So $\gsvec ''\cT$ is well-behaved.  Therefore $\cT^M$ is well-behaved,
a contradiction.
\end{proof}
\end{proof}
\section{An Erd\"os property}

It has been observed that the existence of sharps for all sets 
(rather than a measurable cardinal $\gO$) suffices to
develop the theory in \cite{St1}, and consequently, 
to prove a versions of the results in this paper.
Theorem~\ref{erdos} below, an improvement to Theorem~\ref{wc},
illustrates how an Erd\"os partition property on $\gO$ consistent with 
$V=L$ suffices in many core model applications.

Let $\Ktil$ be the premouse determined by
the inductive definition of the core model in \cite[\S6]{St1}.
By \cite{St1}, if $\gO$ is a measurable cardinal,
then $K = \Ktil$ up to $\gO$, but it is not known if one can show in ZFC
that $\Ktil$ is a proper class.

\begin{theorem}\label{erdos}
Suppose that there is no inner model with a Woodin cardinal.
and that $\gk < \gO$ are cardinals satisfying the partition relation
$\gO \lra ( \go)^{<\go}_\gk$.
Suppose that $\gk$ is a weakly compact cardinal.
Then $\Ktil$ is a premouse of height 
$> \gk^+$ and $(\gk^+ )^{\Ktil} = \gk^+$.
\end{theorem}

Of course, Theorem~\ref{wc} is just a special case of Theorem~\ref{erdos}.
The same method of proof can be used to show that the hypothesis
$$\gO \lra ( \go)^{<\go}_\gk$$ for all $\gk < \gO$,
is more than enough
for all applications in \cite{MiSch}, \cite{MiSchSt}, \cite{Sch1}, 
\cite{Sch2}, \cite{St1}, and \cite{St3}.
In fact, the proof uses only a model theoretic consequence of the 
partition relation.

\begin{proof}[Proof of Theorem~\ref{erdos}]
Suppose that $\gk < \gO$ are cardinals such that
$\gk$ is a weakly compact and
$$\gO \lra ( \go )^{<\go}_\gk$$
Assume that there is no inner model with a Woodin cardinal.

We shall adapt an argument of Reinhardt and Silver;
see \cite[pp.394--5]{J}.
Let $\fA$ be the structure 
$\langle \, V_\gO \, , \in \, , \gk \, \rangle$
expanded to include predicates for Skolem functions.
Let $I$ be an infinite set of indiscernibles for $\fA$
such that $I \subset (\gk , \gO)$.
The existence of $I$ follows from our assumption that 
$\gO \lra (\go )^{<\go}_\gk$ in the standard way.
Let $$\gp: \fB \lra \fA$$ be the inverse transitive collapse of
the closure of $I$ under the Skolem functions for $\fA$.
Fix some $\ga \in I$
and put $\gabar = \gp^{-1}(\ga )$ and $\gkbar = \gp^{-1}(\gk )$.
Let $j: \fB \lra \fB$ 
be an elementary embedding obtained by ``shifting'' $\gabar$,
so that $\crit (j) \leq \gabar$.
With a judicious choice of $I$, 
we can arrange that $\crit (j) > \gkbar$.

Now, using $j$ just as the embedding
from a normal measure was used in \cite{St1},
one can show that
the inductive definition of the core model, as carried out in $\fB$,
succeeds up to $\crit (j)$.
The main 
point is that the superstrong extender from $j$ is weakly amenable
to $\fB$.
Since $\gp$ is an elementary embedding, $\Ktil$ has height 
$\geq j ( \crit (j) ) > \gk^+$.
Arguing as in Theorem~\ref{wc}, we see that
$$\fB \models (  \gkbar \,^+ )^K = \gkbar \,^+ \ .$$
By elementarity again, we conclude that $(\gk^+ )^{\Ktil} = \gk^+$.
\end{proof}

\section{The core model in $L[\Evec ]$}\label{section-LE}
 
This section is still in preparation.

\end{document}